\documentclass[12pt]{article}

\usepackage{amsfonts, latexsym}

\newcommand{\trm}{\textrm}

\newcommand{\Mobius}{M\"{o}bius}
\newcommand{\blst}{\begin{trivlist}}
\newcommand{\elst}{\end{trivlist}}
\newcommand{\bibi}{\bibitem}

\newtheorem{thm}{Theorem}[section]
\newtheorem{prop}[thm]{Proposition}
\newtheorem{cor}[thm]{Corollary}
\newtheorem{lem}[thm]{Lemma}
\newtheorem{conj}[thm]{Conjecture}
\newtheorem{exa}[thm]{Example}
\newtheorem{defn}[thm]{Definition}

\newcommand{\ben}{\begin{enumerate}}
\newcommand{\een}{\end{enumerate}}
\newcommand{\ble}{\begin{lem}}
\newcommand{\ele}{\end{lem}}
\newcommand{\bth}{\begin{thm}}
\newcommand{\eth}{\end{thm}}
\newcommand{\bpr}{\begin{prop}}
\newcommand{\epr}{\end{prop}}
\newcommand{\bco}{\begin{cor}}
\newcommand{\eco}{\end{cor}}
\newcommand{\bcon}{\begin{conj}}
\newcommand{\econ}{\end{conj}}
\newcommand{\bde}{\begin{defn}}
\newcommand{\ede}{\end{defn}}
\newcommand{\bex}{\begin{exa}\rm}
\newcommand{\eex}{\end{exa}}
\newcommand{\barr}{\begin{array}}
\newcommand{\earr}{\end{array}}
\newcommand{\btab}{\begin{tabular}}
\newcommand{\etab}{\end{tabular}}
\newcommand{\beq}{\begin{equation}}
\newcommand{\eeq}{\end{equation}}
\newcommand{\beqq}{$$}
\newcommand{\eeqq}{$$}
\newcommand{\bea}{\begin{eqnarray*}}
\newcommand{\eea}{\end{eqnarray*}}
\newcommand{\beaa}{\begin{eqnarray}}
\newcommand{\eeaa}{\end{eqnarray}}
\newcommand{\bce}{\begin{center}}
\newcommand{\ece}{\end{center}}
\newcommand{\bpi}{\begin{picture}}
\newcommand{\epi}{\end{picture}}
\newcommand{\bfi}{\begin{figure} \begin{center}}
\newcommand{\efi}{\end{center} \end{figure}}
\newcommand{\capt}{\caption}
\newcommand{\bsl}{\begin{slide}{}}
\newcommand{\esl}{\end{slide}}

\newcommand{\ch}{\choose}
\newcommand{\bib}{thebibliography}
\newcommand{\pf}{{\bf Proof.}\quad}
\newcommand{\qed}{\rule{1ex}{1ex}}
\newcommand{\Qed}{\hfill{\rule{1ex}{1ex} \medskip}}

\newcommand{\emp}{\emptyset}

\newcommand{\sbe}{\subseteq}

\newcommand{\spe}{\supseteq}

\newcommand{\iso}{\cong}

\newcommand{\mh}{\hat{0}}
\newcommand{\Mh}{\hat{1}}

\newcommand{\zh}{\hat{0}}

\newcommand{\ptn}{\vdash}

\newcommand{\jn}{\vee}
\newcommand{\Jn}{\bigvee}
\newcommand{\mt}{\wedge}

\newcommand{\ree}[1]{(\ref{#1})}

\newcommand{\ra}{\rightarrow}
\newcommand{\Ra}{\Rightarrow}

\newcommand{\Rla}{\Leftrightarrow}

\newcommand{\al}{\alpha}
\newcommand{\be}{\beta}

\newcommand{\de}{\delta}

\newcommand{\om}{\omega}
\newcommand{\si}{\sigma}

\newcommand{\ze}{\zeta}

\newcommand{\De}{\Delta}

\newcommand{\bs}{{\bf s}}

\newcommand{\bu}{{\bf u}}
\newcommand{\bv}{{\bf v}}
\newcommand{\bw}{{\bf w}}
\newcommand{\bx}{{\bf x}}
\newcommand{\by}{{\bf y}}

\newcommand{\bbN}{{\mathbb N}}

\newcommand{\bbR}{{\mathbb R}}
\newcommand{\bbZ}{{\mathbb Z}}
\newcommand{\cA}{{\cal A}}

\newcommand{\cW}{{\cal W}}

\newcommand{\scl}{\scriptstyle}

\newcommand{\cho}{\choose}

{}
\newcommand{\aim}{Adv.\ in Math.}{}
\newcommand{\alu}{Alg.\ Univ.}{}
{}
{}
\newcommand{\dm}{Discrete Math.}{}
{}
{}
{}
{}
{}
{}
{}
\newcommand{\jct}{J.\ Combin.\ Theory}{}
\newcommand{\jcta}{J.\ Combin.\ Theory Ser.\ A}{}
{}
{}
{}
{}
{}
{}
{}
{}
{}
{}

\usepackage{pspic}

\newcommand{\ps}{1}
\newcommand{\pointsize}[1]{\renewcommand{\ps}{#1}}

\newcommand{\Gda}{\put(30,0){\circle*{\ps}}}

\newcommand{\Gac}{\put(0,20){\circle*{\ps}}}
\newcommand{\Gbc}{\put(10,20){\circle*{\ps}}}
\newcommand{\Gcc}{\put(20,20){\circle*{\ps}}}

\newcommand{\Gec}{\put(40,20){\circle*{\ps}}}

\newcommand{\Ggc}{\put(60,20){\circle*{\ps}}}

\newcommand{\Gae}{\put(0,40){\circle*{\ps}}}

\newcommand{\Gce}{\put(20,40){\circle*{\ps}}}
\newcommand{\Gde}{\put(30,40){\circle*{\ps}}}

\newcommand{\Gfe}{\put(50,40){\circle*{\ps}}}
\newcommand{\Gge}{\put(60,40){\circle*{\ps}}}

\newcommand{\Gdg}{\put(30,60){\circle*{\ps}}}

\newcommand{\GacL}[2]{\Gac \put(0,20){\makebox(0,0){\hspace{#1}#2}}}
\newcommand{\GbcL}[2]{\Gbc \put(10,20){\makebox(0,0){\hspace{#1}#2}}}
\newcommand{\GccL}[2]{\Gcc \put(20,20){\makebox(0,0){\hspace{#1}#2}}}

\newcommand{\GecL}[2]{\Gec \put(40,20){\makebox(0,0){\hspace{#1}#2}}}

\newcommand{\GgcL}[2]{\Ggc \put(60,20){\makebox(0,0){\hspace{#1}#2}}}

\newcommand{\GaeL}[2]{\Gae \put(0,40){\makebox(0,0){\hspace{#1}#2}}}

\newcommand{\GceL}[2]{\Gce \put(20,40){\makebox(0,0){\hspace{#1}#2}}}
\newcommand{\GdeL}[2]{\Gde \put(30,40){\makebox(0,0){\hspace{#1}#2}}}

\newcommand{\GfeL}[2]{\Gfe \put(50,40){\makebox(0,0){\hspace{#1}#2}}}
\newcommand{\GgeL}[2]{\Gge \put(60,40){\makebox(0,0){\hspace{#1}#2}}}

\newcommand{\Gdaac}{\put(30,0){\line(-3,2){30}}}
\newcommand{\Gdabc}{\put(30,0){\line(-1,1){20}}}
\newcommand{\Gdacc}{\put(30,0){\line(-1,2){10}}}

\newcommand{\Gdaec}{\put(30,0){\line(1,2){10}}}

\newcommand{\Gdagc}{\put(30,0){\line(3,2){30}}}

\newcommand{\Gacae}{\put(0,20){\line(0,1){20}}}

\newcommand{\Gacce}{\put(0,20){\line(1,1){20}}}

\newcommand{\Gacfe}{\put(0,20){\line(5,2){50}}}

\newcommand{\Gbcae}{\put(10,20){\line(-1,2){10}}}

\newcommand{\Gbcde}{\put(10,20){\line(1,1){20}}}

\newcommand{\Gbcge}{\put(10,20){\line(5,2){50}}}

\newcommand{\Gccce}{\put(20,20){\line(0,1){20}}}
\newcommand{\Gccde}{\put(20,20){\line(1,2){10}}}

\newcommand{\Gecde}{\put(40,20){\line(-1,2){10}}}

\newcommand{\Gecfe}{\put(40,20){\line(1,2){10}}}

\newcommand{\Ggcfe}{\put(60,20){\line(-1,2){10}}}
\newcommand{\Ggcge}{\put(60,20){\line(0,1){20}}}

\newcommand{\Gaedg}{\put(0,40){\line(3,2){30}}}

\newcommand{\Gcedg}{\put(20,40){\line(1,2){10}}}

\newcommand{\Gdedg}{\put(30,40){\line(0,1){20}}}

\newcommand{\Gfedg}{\put(50,40){\line(-1,1){20}}}

\newcommand{\Ggedg}{\put(60,40){\line(-3,2){30}}}

\begin{document}

\title{Left-modular Elements}
\author{Shu-Chung Liu\\
Institute of Mathematics\\
Academia Sinica, Nankang\\
Taipei, Taiwan, R.O.C.\\
liularry@math.sinica.edu.tw\\
\\
Bruce E. Sagan\\
Department of Mathematics\\
Michigan State University\\
East Lansing, MI 48824 U.S.A.\\
sagan@math.msu.edu\\
\\
\emph{Dedicated to the memory of Gian-Carlo Rota}\\
\emph{without whose work on M\"obius functions}\\
\emph{this paper might never have been written}}

\date{\today\\[1in]
	\begin{flushleft}
	\end{flushleft}
     }
\maketitle

\begin{flushleft} Proposed running head: \end{flushleft}
	\begin{center} 
Left-modular elements
	\end{center}

Send proofs to:
\begin{center}
Bruce E. Sagan \\ Department of Mathematics \\Michigan State
University \\ East Lansing, MI 48824-1027\\[5pt]
Tel.: 517-355-8329\\
FAX: 517-432-1562\\
Email: sagan@math.msu.edu
\end{center}

\pagebreak

  \begin{abstract}
Left-modularity~\cite{BS:mfl97}
is a more general concept than modularity in lattice theory.
In this paper, we give a characterization of left-modular elements and
demonstrate two formulae for the characteristic polynomial of a
lattice with such an
element, one of which generalizes Stanley's Partial Factorization Theorem.
Both formulae provide us with  inductive proofs for
Blass and Sagan's Total Factorization Theorem for LL lattices.
The characteristic polynomials and the \Mobius{} functions of non-crossing
partition lattices and shuffle posets are computed as examples.
  \end{abstract}

\pagebreak

\section{Left-modular elements}

Throughout this paper $L$ is a finite lattice where $\mh=\mh_L$ and 
$\Mh=\Mh_L$ are the
minimal and maximal elements, respectively.
We say that $x$ is {\it covered\/} by $y$, and write $x\prec y$, if $x<y$ and
there is no element $z\in L$ such that $x<z<y$.

We use $\mt$ for the meet (greatest lower bound) and $\jn$ for the join
(least upper bound) in $L$.
Given any $x,y,z\in L$ with $z<y$, the \emph{modular inequality}
\beq
z\jn (x\mt y)\le (z\jn x)\mt y   \label{ie_modp}
\eeq
is always true and equality holds whenever $y$ or $z$ is
comparable to $x$.
We say that $x$ and $y$ form a \emph{modular pair} $(x,y)$ if 
\ree{ie_modp} is an equality for any $z<y$.
Note that this relation is not symmetric, in general.
Two kinds of elements are associated to the modular pair:

\bde
\ben
\item An element $x$ is called a \emph{left-modular element} if
$(x,y)$ is a modular pair for every $y\in L$.
\item An element $x$ is called a \emph{modular element} if
both $(x,y)$ and $(y,x)$ are modular pairs for every $y\in L$.
\een
\ede

In a semimodular lattice with rank function $\rho$, the pair $(x,y)$ is
modular if and only if $\rho(x\mt y)+\rho(x\jn y)=\rho(x)+\rho(y)$
\cite[p.~83]{Bir:lt67}; so in this case the relation of being a modular
pair is symmetric, and then there is no difference between modularity and
left-modularity.
However, there are examples such as the non-crossing partition lattices
(see Sec.~\ref{s_ncp}) and the Tamari lattices where the two concepts do
not coincide.

Let $L$ be a graded lattice of rank $n$ with rank function $\rho$.
Then the \emph{characteristic polynomial} of $L$ is defined by
$$
\chi(L,t) =\sum_{x\in L} \mu(x)t^{n-\rho(x)}
$$
where $t$ is an indeterminate, $\mu:L\times L\ra\bbZ$ 
is the M\"obius function of $L$, and $\mu(x)=\mu(\zh,x)$.
There are two important factorization theorems for $\chi$ given by
R.~Stanley:

\bth[Partial Factorization Theorem~\cite{Sta:meg71}] \label{t_mfac}
Let $L$ be an atomic, semimodular lattice (i.e., a geometric lattice)
of rank $n$.
If $x$ is a modular element of $L$, then
\beqq
\chi(L,t)= \chi([\mh,x],t) \sum_{b\ :\ b\mt x=\mh}
\mu(b) t^{n-\rho(x)-\rho(b)}.\ \Qed
\eeqq
\eth
\bth[Total Factorization Theorem~\cite{Sta:sl72}]  \label{t_mlfac}
Let $(L,\De)$ be a supersolvable, semimodular lattice of rank $n$ with
$\De:\mh=x_0\prec x_1\prec\ldots\prec x_n=\Mh$.
Then
\beq        \label{e_mlfac}
\chi(L,t)=(t-a_1)(t-a_2)\cdots (t-a_n)
\eeq
where $a_i$ is the number of atoms of $L$ that are below $x_i$ but
not below $x_{i-1}$.\ \qed
\eth

Note that all elements in the maximal chain $\De$ of a supersolvable lattice
are left-modular (see \cite{Sta:meg71}).
So the hypotheses of Theorem~\ref{t_mlfac} imply that they are modular.
In recent work~\cite{BS:mfl97}, A.~Blass and B.~Sagan generalized the
Total Factorization Theorem to LL lattices where the first ``L'' stands
for the fact that the lattice has a maximal chain all of whose elements are
all left-modular. 
The purpose of this paper is to generalize the Partial Factorization Theorem
by replacing the modular element with a left-modular one and relaxing
the hypotheses requiring that the lattice be atomic and semimodular.
To do so, we will derive a general characterization of  left-modular
elements in this section.
In the next section, we introduce a generalized rank function for a
lattice which might not be graded in the usual sense, and then develop
a general formula for the 
characteristic polynomial of a lattice with a left-modular element in
Theorem~\ref{t_lmgr}.
Under an extra rank-preserving hypothesis we obtain our generalization
of the Partial Factorization Theorem (Theorem~\ref{t_lmfac}).
In Sections~\ref{s_ncp} and \ref{s_shp}, we calculate the characteristic
polynomials and the \Mobius{} functions of the non-crossing partition
lattices and the shuffle posets by using these two formulae, respectively.
The last section contains
two inductive proofs for Blass and Sagan's Total Factorization Theorem for
LL
lattices using our two main theorems.
Consequently, our factorization theorem generalizes the three others.

We say that $y$ is a \emph{complement} of $x$ if $x\mt y=\mh$ and $x\jn
y=\Mh$.
Stanley \cite{Sta:meg71} showed that, in an atomic and semimodular lattice,
$x$ is modular if and only if no two complements of $x$ are comparable.
The next theorem provides an analog for left-modular elements.

\bth        \label{t_comp}
Let $x$ be an element of any lattice $L$.
The following statements are equivalent:
\ben
\item[i.] The element $x$ is left-modular.
\item[ii.] For any $y$, $z\in L$ with $z<y$, we have
   $x\mt z \ne x\mt y$ or $x\jn z \ne x\jn y.$
\item[iii.] For any $y$, $z\in L$ with $z\prec y $, we have
   $x\mt z= x\mt y$ or $x\jn z= x\jn y$ but not both.
\item[iv.] For every interval $[a,b]$ containing $x$,
   no two complements of $x$ with respect to
   the sublattice $[a,b]$ are comparable.
\een
\eth
\pf
We will prove the implications (i) $\Ra$ (ii) $\Ra$ (iii) $\Ra$ (i).
The proof of (ii) $\Rla$ (iv) is immediate.

First we make some preliminary observations.
Suppose $z<y$.
We claim that $x\jn y=x\jn z$ if and only if $y=(z\jn x)\mt y$.
The forward direction is trivial since $(x\jn y)\mt y=y$.
For the reverse, note that  $y=(z\jn x)\mt y$ implies $y\le x\jn z$.
Now $z<y\le x\jn z$, and joining all sides with $x$ gives $x\jn y=x\jn z$.
Dually $x\mt y=x\mt z$ if and only if $z= z\jn(x\mt y)$.

For any $z<y$ the inequalities
  \beq
z\le z\jn(x\mt y)\le (z\jn x)\mt y\le y   \label{ineq}
  \eeq
are true by the modular inequality~\ree{ie_modp}.
Since $z\ne y$, at least one of the $\le$'s in (\ref{ineq}) should be $<$.
Therefore (i) $\Ra$ (ii).
If $z\prec y$, then exactly two of the $\le$'s should be $=$ and the
remaining one must be $\prec$.
Thus (ii) $\Ra$ (iii).

To show (iii) $\Ra$ (i), let us consider the contrapositive:
assume that there are $u$, $v\in L$ with $u<v$ such that $u\jn (x\mt v)<
(u\jn x)\mt v$.
Given any $y$, $z\in [u\jn (x\mt v),(u\jn x)\mt v]$ with $z\prec y$, we
have $y\le (u\jn x)\mt v\le v$. This implies 
$u\jn (x\mt y)\le u\jn (x\mt v)\le z$,
so that $x\mt y\le z$.
It follows that $x\mt z=x\mt y$.
Similarly, we can get $x\jn z=x\jn y$.
\Qed

The existence of a left-modular element in $L$ implies that such
elements are also present in certain sublattices as the next
proposition shows. 

\bpr        \label{p_lmjm}
Let $x$ be a left-modular element in lattice $L$.
Then for any $y\in L$
\ben
\item the meet $x\mt y$ is a left-modular element in $[\mh,y]$, and
\item the join $x\jn y$ is a left-modular element in $[y,\Mh]$.
\een
\epr
\pf
Let $a$, $b\in [\mh, y]$ with $b<a$.
By left-modularity of $x$, we have
  \bea
b\jn ((x\mt y)\mt a) & = & b\jn (x\mt (y\mt a))\
  =\ (b\jn x)\mt (y\mt a)\\
& = & ((b\jn x)\mt y)\mt a\ =\ (b\jn (x\mt y))\mt a.
  \eea
So $x\mt y$ is a left-modular element in $[\mh,y]$.
The proof for join is similar.\ \Qed

\section{The characteristic polynomial}

We begin with a general lemma.

\ble        \label{l_lmgr}
Let $L$ be a lattice with an arbitrary function $r:L\ra \bbR$ and
let {$n\in \bbR$}.
If $x\in L$ is a left-modular element, then
  \beqq
\sum_{y\in L}\mu(y)t^{n-r(y)}=
\sum_{b\mt x=\mh}\mu(b) \sum_{y\in [b,b\jn x]}
  \mu(b,y) t^{n-r(y)}.   \label{sumy}
  \eeqq
\ele
\pf
We will mimic Stanley's proof in \cite{Sta:meg71}.
By Crapo's Complementation Theorem~\cite{Cra:mfl66}, for any given
$a\in [\mh,y]$
  $$\mu(y)=\sum_{a',a''} \mu(\mh,a') \ze(a',a'') \mu(a'',y),$$
where $a'$ and $a''$ are complements of $a$ in $[\mh,y]$, and $\ze$ is
the zeta function defined by $\ze(u,v)=1$ if $u\le v$ and
$\ze(u,v)=0$ otherwise.
Let us choose $a=x\mt y$.
The element $a$ is left-modular in $[\mh,y]$ by Proposition~\ref{p_lmjm}.
But no two complements of $a$ in $[\mh,y]$ are comparable by
Theorem~\ref{t_comp}.
Thus
  \beq
\mu(y)=\sum_b \mu(\mh,b) \mu(b,y),   \label{muy}
  \eeq
where the sum is over all complements $b$ of $a$ in $[\mh,y]$, i.e.,
over all $b$ satisfying $b\le y$, $b\mt (x\mt y)=\mh$ and
$b\jn (x\mt y)=y$.
Since $x$ is left-modular, it is equivalent to say that
the sum in~(\ref{muy}) is over all $b\in L$ satisfying $b\mt x=\mh$ and
$y\in [b,b\jn x]$.
Thus we have
  \bea
\sum_{y\in L} \mu(y) t^{n-r(y)}
& = & \sum_{y\in L}\quad 
\sum_{\scriptstyle b\mt x=\mh \atop\rule{0pt}{9pt}\scriptstyle y\in [b,b\jn x]}
  \mu(\mh,b) \mu(b,y) t^{n-r(y)}\\
& = & \sum_{b\mt x=\mh} \mbox{\ } \mu(b) \sum_{y\in [b,b\jn x]}
  \mu(b,y) t^{n-r(y)}.\ \Qed
  \eea

Obviously, the previous lemma is true for the ordinary rank function if $L$
is graded.
To apply this result to more general lattices we make the following
definition.

\bde         \label{d_gr}
A \emph{generalized rank function} of a lattice $L$ is a function
$\rho:\{(x,y)\in L\times L\mid x\le y \}\ra \bbR$ such that
for any $a\le b\le c$
  $$\rho(a,c)=\rho(a,b)+\rho(b,c).$$
In this case, we say $L$ is \emph{generalized graded} by $\rho$.
\ede

For short we write $\rho(x)=\rho(\mh,x)$.
Conversely, if we take any function $\rho:L\ra \bbR$ such that
$\rho(\mh)=0$,
then we can easily construct a generalized rank function, namely
$\rho(x,y)=\rho(y)-\rho(x)$.
So the ordinary rank function is a special case.

If $L$ is generalized graded by $\rho$, we now define a generalized
characteristic polynomial of $L$ by
\beq        \label{d_ndefchi}
\chi(L,t)=\sum_{x\in L}\mu(x)t^{\rho(x,\Mh)}=\sum_{x\in L}
\mu(x)t^{\rho(\Mh)-\rho(x)}.
\eeq
Note that $\chi$ will depend on which generalized rank function we pick.
Since the restriction of a generalized rank function to an interval $[a,b]$
still satisfies Definition~\ref{d_gr} with $L=[a,b]$, the characteristic
polynomial of the interval is defined in the same manner.

The following theorem, which follows easily from Lemma~\ref{l_lmgr},
is one of our main results.
In it, the \emph{support} of $\mu$ is defined by
$$
H(L)=\{x\in L\mid \mu(x)\ne 0\}.
$$

\bth        \label{t_lmgr}
Let $L$ be generalized graded by $\rho$.
If $x\in L$ is a left-modular element, then
\beq       \label{e_lmgr}
\chi(L,t)= \sum_{\scl b\in H(L)\atop\rule{0pt}{9pt}\scl b\mt x=\mh}
\left[ \mu(b)t^{\rho(\Mh)-\rho(b\jn x)} \chi([b,b\jn x],t)\right].\ \qed
\eeq
\eth

In the sum~(\ref{e_lmgr}), the term $\chi([b,b\jn x],t)$ depends on $b$.
To get a factorization formula, we will remove the dependency by applying
certain restrictions so that $\chi([b,b\jn x],t)=\chi([\mh,x],t)$ for all
$b$ in the sum.

First, we will obtain a general condition under which two lattices have
the same characteristic polynomial.
In the following discussion, let $L$ and $L'$ be lattices and let
$\tau:L\ra L'$ be any map.
For convenience, we also denote $\mh=\mh_L$, $\mh'=\mh_{L'}$ and similarly
for $\Mh$, $\Mh'$, $\mu$, $\mu'$, etc.

We say $\tau$ is a \emph{join-preserving} map if
$$
\tau(u\jn v)=\tau(u)\jn \tau(v)
$$
for any $u$, $v\in L$.
Note that from this definition $\tau$ is also order-preserving since
$$
x\le y\ \Ra\ y=x\jn y\ \Ra\ \tau(y)=\tau(x\jn y)=\tau(x)\jn\tau(y)\ \Ra\
\tau(x)\le\tau(y).
$$

If $\tau$ is join-preserving, then given any $x'\in\tau(L)$, we claim that
the subset $\tau^{-1}(x')$ has a unique maximal element in $L$.
Suppose that $\tau(u)=\tau(v)=x'$ for some $u$, $v\in L$.
We have $\tau(u\jn v)=\tau(u)\jn\tau(v)=x'$.
Thus $u\jn v\in\tau^{-1}(x')$ and the claim follows.

If, in addition, $\tau$ is surjective then we can define a map
$\si:L'\ra L$ by
\beq        \label{d_simax}
\si(x')=\trm{ the maximal element of $\tau^{-1}(x')$}.
\eeq
The map $\si$ must also be order preserving.  To see this, suppose
$x'\le y'$ in $L'$ and consider $x=\si(x'), y=\si(y')$.  Then
$$
\tau(x\jn y)=\tau(x)\jn\tau(y)=x'\jn y'=y'.
$$
So $x\jn y\in\tau^{-1}(y')$ which forces $x\jn y\le y$ by definition
of $\si$.  Thus $x\le y$ as desired.

\ble       \label{l_musame}
Using the previous notation, suppose that $\tau$ is surjective and
join-preserving and that $\si$ satisfies $\si(\mh')=\mh$.
Then for any $x'\in L'$ we have
  $$\mu'(x')=\sum_{y\in \tau^{-1}(x')}\mu(y).$$
\ele
\pf
This is trivial when $x'=\mh'$.
Let $x=\si(x')$.
{From} the assumptions on $\tau$ and $\si$ it is easy to see that
\beq      \label{e_Tinv}
[\mh,x]= \biguplus_{y'\in[\mh',x']}\tau^{-1}(y').
\eeq
Now, by surjectivity of $\tau$ and induction, we get
$$
\mu'(x')=-\sum_{y'<x'}\mu'(y')
  =-\sum_{\scl y\in\tau^{-1}(y')\atop\rule{0pt}{9pt} \scl y'<x'}\mu(y)
  = \sum_{y\in \tau^{-1}(x')}\mu(y).\ \Qed
$$

Let $L$ and $L'$ be generalized graded by $\rho$ and $\rho'$, respectively.
We say an order-preserving map $\tau:L\ra L'$ is \emph{rank-preserving} on
a subset $S\sbe L$ if $\rho(x,y)= \rho'(\tau(x),\tau(y))$ for any
$x$, $y\in S$, $x\le y$.

\ble       \label{l_chisame}
If, in addition to the hypotheses of Lemma~\ref{l_musame}, the map $\tau$
is rank-preserving on $H(L)\cup\{\Mh\}$ then
$$
\chi(L,t)= \chi(L',t).
$$
\ele
\pf
{From} (\ref{e_Tinv}) in the proof of Lemma~\ref{l_musame}, we know $L=
\biguplus_{x'\in L'}\tau^{-1}(x')$.
Then by Lemma~\ref{l_musame} and the rank-preserving nature of $\tau$,
we have
\bea
\chi(L',t)&=&\sum_{x'\in L'}\mu'(x')t^{\rho'(x',\Mh')}\\
  &=&\sum_{x'\in L'}\ \sum_{y\in \tau^{-1}(x')}
\mu(y)t^{\rho'(x',\Mh')}\\
  &=&\sum_{x'\in L'}\ \, \sum_{y\in \tau^{-1}(x')\cap H(L)}
\mu(y)t^{\rho'(\tau(y),\tau(\Mh))}\\
  &=&\sum_{y\in H(L)}\mu(y)t^{\rho(y,\Mh)}\\
  &=& \chi(L,t).\ \Qed
\eea

It is easy to generalize the previous lemma to arbitrary posets as long as
the map $\si$ is well defined.  However, we know of no application of the
result in this level of generality.

Returning to our factorization theorem, we still need one more tool.
For any given $a$, $b$ in a lattice, we define
$$\si_a :[b,a\jn b]\ra [a\mt b,a]\quad\trm{by}\quad\si_a(u)=u\mt a,$$
$$\tau_b:[a\mt b,a]\ra [b,a\jn b]\quad\trm{by}\quad\tau_b(v)=v\jn b.$$
The map $\tau_b$ is the one we need to achieve
$\chi([b,b\jn x],t)=\chi([\mh,x],t)$.
In the following, we write $H(x,y)$ for $H([x,y])$ which is the support
of $\mu$ defined on the sublattice $[x,y]$.
We can now prove our second main result.

\bth        \label{t_lmfac}
Let $L$ be generalized graded by $\rho$ and
let $x\in L$ be an left-modular element.
If the map $\tau_b$ is rank-preserving on $H(\mh,x)\cup\{x\}$
for every $b\in H(L)$ satisfying $b\mt x=\mh$.
Then
  \beq        \label{e_lmfac}
\chi(L,t)
=\chi([\mh,x],t)\sum_{\scl b\in H(L)\atop\rule{0pt}{9pt} \scl b\mt x=\mh}
   \mu(b)t^{\rho(\Mh)-\rho(x)-\rho(b)}.
  \eeq
\eth
\pf
First, we will show that $\chi([b,b\jn x],t)= \chi([\mh,x],t)$ for any
$b\in H(L)$ with $b\mt x=\mh$ by
verifying the hypotheses of Lemma~\ref{l_chisame}.
By left-modularity of $x$, we have
\beq    \label{e_ts}
\tau_b\si_x(y)= b\jn (x\mt y)=(b\jn x)\mt y=y
\eeq
for any $y\in [b,b\jn x]$.
So $\tau_b$ is surjective.
And it is easy to check that $\tau_b$ is join-preserving.
As for $\si_x$, we must check that it satisfies the
definition~\ree{d_simax}.
Given $z\in\tau_b^{-1}(y)$ we have $y=\tau_b(z)=z\jn b$.
So by the modular inequality~(\ref{ie_modp}) we get
$$
\si_x(y)=y\mt x=(z\jn b)\mt x\ge z\jn (b\mt x)\ge z.
$$
Since this is true for any such $z$, we have $\si_x(y)\ge\max
\tau_b^{-1}(y)$.
But equation~\ree{e_ts} implies $\si_x(y)\in\tau_b^{-1}(y)$, so we
have equality.
Finally $\mh_{[b,b\jn x]}=b$ so $\si_x(b)=b\mt x=\mh$ as desired.

Now we need only worry about the exponent on $t$ in Theorem~\ref{t_lmgr}.
But since $\tau_b$ is rank-preserving on $H(\mh,x)\cup\{x\}$, we get
$$
\rho(b\jn x) = \rho(\mh,b)+\rho(b,b\jn x)
   = \rho(\mh,b)+\rho(\mh,x)= \rho(b)+\rho(x).\ \Qed
$$

Here we state a corollary which relaxes the hypothesis in Stanley's Partial
Factorization Theorem.

\bco      \label{c_mfac}
Equation~\ree{e_lmfac} holds when $L$ is a semimodular lattice (graded by
the ordinary rank function) with a modular element $x$.
\eco
\pf
To apply Theorem~\ref{t_lmfac},
it suffices to show that $\rho(\mh,z)=\rho(b,z\jn b)$ for every
$z\in [\mh,x]$.
Since $(b,x)$ is a modular pair, we have $(z\jn b)\mt x=z\jn (b\mt x)
=z\jn\mh=z$.
By Proposition~\ref{p_lmjm}, $z=(z\jn b)\mt x$ is left-modular in
$[\mh,z\jn b]$, so $(z,b)$ is a modular pair in this lattice.
Thus $\rho(z\mt b)+\rho(z\jn b)=\rho(z)+\rho(b)$, because
$[\mh,z\jn b]$ is a semimodular lattice.
Since $z\mt b=\mh$ we are done.
\Qed

We take the divisor lattice $D_n$ as an example.
It is semimodular, but not atomic in general, so Stanley's theorem does
not apply.
However, Corollary~\ref{c_mfac} can be used for any $x\in D_n$, since all
elements are modular.

We will now present a couple of applications of the previous results in the
following two sections.

\section{Non-crossing Partition Lattices}   \label{s_ncp}

The non-crossing partition lattice was first studied by
Kreweras~\cite{Kre:pnc72} who showed its \Mobius{} function is related to
the Catalan numbers.
By using NBB sets (see Sec.~\ref{s_lft} for the definition), Blass and
Sagan~\cite{BS:mfl97} combinatorially explained this fact.
In this section we will calculate the characteristic polynomial for
a non-crossing partition lattice and then offer another explanation for
the value of its \Mobius{} function.

If it causes no confusion, we will not explicitly write out any blocks
of a partition that are singletons.
Let $n\ge 1$.
We say that a partition $\pi\ptn [n]$ is \emph{non-crossing} if there do
not exist two distinct blocks $B, C$ of $\pi$ with
$i$, $k\in B$ and $j$, $l\in C$ such that $i<j<k<l$.
Otherwise $\pi$ is \emph{crossing}.

Another way to view non-crossing partitions will be useful.
Let $G=(V,E)$ be a graph with vertex set $V=[n]$ and edge set $E$.
We say that $G$ is \emph{non-crossing} if, when the vertices are arranged in
their natural order clockwise around a circle and the edges are drawn as
straight line segments, no two edges of $G$ cross geometrically.
Given a partition $\pi$ we can form a graph $G_{\pi}$ by representing each
block $B=\{i_1<i_2<\ldots<i_l\}$ by a cycle with edges
$i_1i_2,i_2i_3,\ldots,i_li_1$.
(If $|B|=1$ or 2 then $B$ is represented by an isolated vertex or edge,
respectively.)
Then it is easy to see that $\pi$ is non-crossing as a partition if and
only if $G_\pi$ is non-crossing as a graph.

The set of non-crossing partitions of $[n]$, denoted by $NC_n$, forms a
meet-sublattice
of partition lattice $\Pi_n$ with the same rank function.
However unlike $\Pi_n$, the non-crossing partition lattice is not
semimodular in general, since if
$\pi=13$ and $\si=24$ then $\pi\mt\si=\mh$ and $\pi\jn\si=1234$.
So we have
  $$\rho(\pi)+\rho(\si)=2<3=\rho(\pi\mt\si)+\rho(\pi\jn\si).$$
The $\Pi_n$-join $\pi \jn \si=13/24$ also explains why
$NC_n$ is not a sublattice of $\Pi_n$.

Let $n\ge 2$ and $\pi=12\ldots (n-1)$.
It is well-known~\cite{Sta:sl72} that $\pi$ is modular in $\Pi_n$ and so
left-modular there.
Given any $\al$, $\be\in NC_n$ with $\al<\be$ and both incomparable to $\pi$.
It is clear that $\al\jn\pi=\be\jn\pi=\Mh$ in $\Pi_n$ as well as in $NC_n$.
By Theorem~\ref{t_comp} we get $\al\mt\pi<\be\mt\pi$ in $\Pi_n$.
Since $NC_n$ is a meet-sublattice of $\Pi_n$, this inequality for
the two meets still holds in $NC_n$.
This fact implies that $\pi$ is left-modular in $NC_n$.
In general, $\pi$ is not modular in $NC_n$.
If $n\ge 4$, let $\si =2n$ and $\phi =1(n-1)/23\ldots (n-2)$.
Clearly $\phi <\pi$, $\pi\mt\si=\phi\mt\si=\mh$ and
$\pi\jn\si=\phi\jn\si=\Mh$ in $NC_n$, so that $(\si,\pi)$ is not a modular
pair.

\bpr        \label{p_cpnc}
The characteristic polynomial of the non-crossing partition lattice
$NC_n$ satisfies
  $$\chi(NC_n,t)=t\ \chi(NC_{n-1},t)-
\sum_{i=1}^{n-1}\chi(NC_i,t)\chi(NC_{n-i},t)$$
with the initial condition $\chi(NC_1,t)=1$.
\epr
\pf
The initial condition is trivial.
Let $n\ge 2$ and $\pi=12\ldots (n-1)$.
We will apply Theorem~\ref{t_lmgr}.
Note that $b\mt\pi=\mh$ if and only if any two numbers of $[n-1]$ are in
different blocks of $b$, so either $b=\mh$ or $b=mn$ with $1\le m\le n-1$.

If $b=\mh$, then $\chi([b,b\jn\pi],t)=\chi([\mh,\pi],t)=\chi(NC_{n-1},t)$.
Thus we get the first term of the formula.
Now let $b=mn$.
It is clear that $b\jn\pi=\Mh$, so we need to consider the sublattice
$[b,\Mh]$.
Given any $\om\in [b,\Mh]$, the edge $mn$ (which may not be in $E(G_{\om})$)
geometrically separates the graph $G_{\om}$ into two parts, $G_{\om,1}$
and $G_{\om,2}$, which are induced by vertex sets $\{1,2,\ldots,m,n\}$
and $\{m,m+1,\ldots,n-1,n\}$, respectively.
By contracting the vertices $m$ and $n$ in both $G_{\om,1}$ and $G_{\om,2}$,
we get two non-crossing graphs $\bar{G}_{\om,1}$ and $\bar{G}_{\om,2}$.
It is easy to check that the map $f:[b,\Mh]\ra NC_m\times NC_{n-m}$
defined by $f(G_{\om})=(\bar{G}_{\om,1},\bar{G}_{\om,2})$ is an isomorphism
between these two lattices.
Therefore
  $$\chi([b,b\jn\pi],t)=\chi(NC_m,t)\chi(NC_{n-m},t),$$
and the proof is complete.
\ \Qed

For any $\om=B_1/B_2/\ldots/B_k\in NC_n$, the interval
$[\mh,\om]\iso\prod_i NC_{|B_i|}$.
Hence to compute the \Mobius{} function of $NC_n$, it suffices to do this
only for $\Mh$.
By Proposition~\ref{p_cpnc} we have the recurrence relation
  \bea
\mu(NC_n)&=&\chi(NC_n,0)\\
&=&-\sum_{i=1}^{n-1}\chi(NC_i,0)\chi(NC_{n-i},0)\\
&=&-\sum_{i=1}^{n-1}\mu(NC_i)\mu(NC_{n-i})
  \eea
with the initial condition $\mu(NC_1)=1$.
Recall that the  Catalan numbers $C_n=\frac{1}{n+1} {2n\ch n}$
satisfy the recurrence relation
  $$C_n= \sum_{i=0}^{n-1}C_i C_{n-1-i}$$
with the initial condition $C_0=1$.
Therefore, by induction, we obtain Kreweras' result that
  $$\mu(NC_n)=(-1)^{n-1}C_{n-1}.$$


\section{Shuffle Posets}     \label{s_shp}

The poset of shuffles was introduced by Greene~\cite{Gre:ps88},
and he obtained a formula for its characteristic polynomial
$$
\chi(\cW_{m,n},t)=(t-1)^{m+n}\sum_{i\ge 0}{m \ch i}{n \ch
i}\frac{1}{(1-t)^i}.
$$
In this section we will derive an equivalent formula by using
Theorem~\ref{t_lmfac}.
Before doing this, we need to recall some definitions and results of Greene.
Let $\cA$ be a set, called the \emph{alphabet of letters}.
A \emph{word} over $\cA$ is a sequence $\bu=u_1 u_2\ldots u_n$ of distinct
letters of $\cA$.
We will sometimes also use $\bu$ to stand for the set of letters in the
word,
depending upon the context.
A \emph{subword} of $\bu$ is $\bw=u_{i_1}\ldots u_{i_l}$ where
$i_1<\ldots<i_l$.
If $\bu$, $\bv$ are any two words then the \emph{restriction} of
$\bu$ to $\bv$ is the subword $\bu_{\bv}$ of $\bu$ whose letters are
exactly those of $\bu\cap\bv$.
A \emph{shuffle} of $\bu$ and $\bv$ is any word $\bs$ such that
$\bs=\bu\uplus\bv$ as sets and $\bs_{\bu}=\bu$,
$\bs_{\bv}=\bv$ as words.

Given nonnegative integers $m$ and $n$, fix disjoint words
$\bx=x_1\ldots x_m$ and $\by=y_1\ldots y_n$.
The \emph{poset of shuffles} $\cW_{m,n}$ consists all shuffles
$\bw$ of a subword of $\bx$ with a subword of $\by$ while
the partial order is that $\bv\le\bw$ if $\bv_{\bx}\spe\bw_{\bx}$,
$\bv_{\by}\sbe\bw_{\by}$ as sets and $\bv_{\bw}=\bw_{\bv}$ as words.
The covering relation is more intuitive:
$\bv\prec\bw$ if $\bw$ can be obtained from $\bv$ by either adding a single
$y_i$ or deleting a single $x_j$.
It is easy to see that $\cW_{m,n}$ has $\mh=\bx$, 
$\Mh=\by$, and is graded by the rank function
$$\rho(\bw)=(m-|\bw_{\bx}|)+|\bw_{\by}|.$$
For example, $\cW_{2,1}$ is shown in Figure~\ref{W21}
where $\bx=de$ and $\by=D$.

\thicklines
\setlength{\unitlength}{2.4pt}
\pointsize{1}
\bfi
\bpi(60,70)(0,-10)
\Gda
\put(30,-5){\makebox(0,0){$de$}}
\GacL{-15pt}{$d$}
\GbcL{-15pt}{$e$}
\GccL{30pt}{$Dde$}
\GecL{30pt}{$dDe$}
\GgcL{30pt}{$deD$}
\GaeL{-15pt}{$\emp$}
\GceL{-25pt}{$Dd$}
\GdeL{25pt}{$De$}
\GfeL{-25pt}{$dD$}
\GgeL{25pt}{$eD$}
\Gdg
\put(30,65){\makebox(0,0){$D$}}
\Gdaac \Gdabc \Gdacc \Gdaec \Gdagc
\Gacae \Gacce \Gacfe \Gbcae \Gbcde \Gbcge
\Gccce \Gccde \Gecde \Gecfe \Ggcfe \Ggcge
\Gaedg \Gcedg \Gdedg \Gfedg \Ggedg
\epi
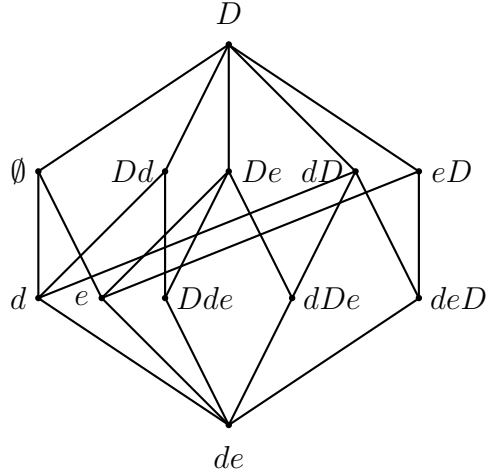
\capt{The lattice $\cW_{2,1}$}     \label{W21}
\efi

Every shuffle poset is actually a lattice.
To describe the join operation in $\cW_{m,n}$, Greene defined crossed
letters as follows.
Given $\bu$, $\bv\in\cW_{m,n}$ then $x\in\bu\cap\bv\cap\bx$ is
\emph{crossed} in $\bu$ and $\bv$ if there exist letters $y_i$,
$y_j\in\by$ with
$i\leq j$ and $x$ appears before $y_i$ in one of the two words but after
$y_j$ in the other.
For example, let $\bx=def$ and $\by=DEF$.
Then in the two shuffles $\bu=dDEe$, $\bv=Fdef$, the only crossed
letter is $d$.
The join of $\bu$, $\bv$ is then the unique word $\bw$ greater than both
$\bu$, $\bv$ such that
  \beqq        \label{Wjn}
\barr{rcl}
\bw_{\bx}&=&\{x\in\bu_{\bx}\cap\bv_{\bx} \mid
  x \trm{ is not crossed}\}\\
\bw_{\by}&=&\bu_{\by}\cup\bv_{\by}.
\earr
  \eeqq
In the previous example, $\bu\jn\bv=DEFe$.
This join also shows that $\cW_{m,n}$ is not semimodular in general,
because
  $\rho(\bu)+\rho(\bv)=3+1<5=\rho(\bu\jn\bv)\le
\rho(\bu\jn\bv)+\rho(\bu\mt\bv).$
Since $(\cW_{n,m})^*=\cW_{m,n}$, the meet operation in $\cW_{m,n}$ is as
same as the join operation in $(\cW_{n,m})^*$.
So to find the meet in the analogous way we need to consider those letter
$y\in\bu\cap\bv\cap\by$ crossed in $\bu$ and $\bv$.

Greene also showed that subwords of $\bx$ and subwords of $\by$ are
modular elements of $\cW_{m,n}$.
In particular, the empty set $\emp$ is modular.
Also note that $[\mh,\emp]\iso B_m$.
We now give our formula for the characteristic polynomial of $\cW_{m,n}$.

\bpr        \label{p_cpsp}
The characteristic polynomial of the shuffle poset is
  \beq        \label{e_cpsp}
\chi(\cW_{m,n},t)=(t-1)^m \sum_{i=0}^n (-1)^i
  {n \ch i}{m+i \cho i} t^{n-i}.
  \eeq
\epr
\pf
Consider any $\bu$ with $\bu\mt\emp=\mh$.
In general, if $\bu\mt\emp=\bw$ then
$\bw_{\bx}=\bu_{\bx}\cup\emp_{\bx}=\bu_{\bx}$.
So $\bu\mt\emp=\mh$ if and only if $\bx$ is a subword of $\bu$,
i.e., the element $\bu$ is a shuffle of $\bx$ with a subword of $\by$.
Furthermore, for any $\bv\in [\mh,\emp]$, there is no crossed letter $x$
in $\bu$ and $\bv$ since $\bv_{\by}=\emp$.
It follows that $(\bu\jn\bv)_{\bx}=\bu_{\bx}\cap\bv_{\bx}=\bv$ and
$(\bu\jn\bv)_{\by}=\bu_{\by}\cup\bv_{\by}=\bu_{\by}$ as sets.
Then we get
  \bea
\rho(\bu\jn\bv)-\rho(\bu)&=&[(m-|\bv|)+|\bu_{\by}|]-[(m-m)+|\bu_{\by}|]\\
&=& m-|\bv|\ =\ \rho(\bv)-\rho(\mh).
  \eea
Thus the map $\tau_{\bu}:[\mh,\emp]\ra [\bu,\emp\jn\bu]$ is rank-preserving.

Since $[\mh,\emp]\iso B_m$ we get, by Theorem~\ref{t_lmfac},
  $$\chi({\cW}_{m,n},t)=(t-1)^m \sum_{\bu\mt\emp=\mh}
  \mu(\bu)t^{(m+n)-\rho(\bu)-m}.$$
It is easy to see that the interval $[\mh,\bu]$ is isomorphic to $B_i$
where $i=|\bu_{\by}|$.
So $\mu(\bu)=(-1)^{|\bu_{\by}|}=(-1)^{\rho(\bu)}$.
Now we conclude that
  \bea
\chi({\cW}_{m,n},t)
&=&(t-1)^m \sum_{i=0}^n\left[{\trm{the number of ways to} \atop
  \trm{shuffle $\bx$ with $i$ letters of $\by$}} \right]
  (-1)^i t^{n-i}\\
&=&(t-1)^m \sum_{i=0}^n (-1)^i{n \ch i}{m+i \cho i} t^{n-i}.\ \Qed
  \eea

To determine the \Mobius{} function of $\cW_{m,n}$, it suffices to compute
$\mu(\Mh)$ since for any $\bw\in\cW_{m,n}$ the interval $[\mh,\bw]$ is
isomorphic to a product of $\cW_{p,q}$'s for certain $p\le m$ and $q\le n$.
Simply plugging $t=0$ into formula~(\ref{e_cpsp}) gives us the \Mobius{}
function $\mu(\cW_{m,n})$.

\bco[Greene, \cite{Gre:ps88}]
We have
  $$\mu(\cW_{m,n})=(-1)^{m+n}{m+n \cho n}.\ \qed$$
\eco


\section{NBB Sets and Factorization Theorems}  \label{s_lft}

Blass and Sagan~\cite{BS:mfl97} derived a Total Factorization Theorem for
LL lattices which generalizes Theorem~\ref{t_mlfac}.
Applying Theorem~\ref{t_lmgr} and \ref{t_lmfac}, respectively, we will
offer two inductive proofs for their theorem.
First of all, we would like to outline their work.

Given a lattice $L$, let $A=A(L)$ is the set of atoms of $L$.
Let $\unlhd$ be an arbitrary partial order on $A$.
A nonempty set $D\sbe A$ is \emph{bounded below} or \emph{BB} if, for
every $d\in D$ there is an $a\in A$ such that
$$a \lhd d \qquad\trm{and}\qquad a < \Jn D.$$
A set $B\sbe A$ is called \emph{NBB} (\emph{no bounded below} subset) if it
does not contain any $D$ which is bounded below.
An NBB set is said to be a base for its join.
One of the main results of Blass and Sagan's paper is the following
theorem which is a simultaneous generalization of both Rota's NBC and
Crosscut Theorems (for the crosscut $A(L)$).

\bth[Blass and Sagan, \cite{BS:mfl97}]   \label{t_NBBmu}
Let $L$ be a finite lattice and let $\unlhd$ be any partial order on $A$.
Then for all $x\in L$ we have
$$\mu(x)= \sum_{B} (-1)^{|B|}$$
where the sum is over all NBB bases $B$ of $x$.
\ \qed
\eth

Given an arbitrary lattice $L$,
let $\De: \mh=x_0\prec x_1\prec \ldots\prec x_n=\Mh$ be a maximal chain of
$L$.
The $i^{\mathit{th}}$ \emph{level} of $A$ is defined by
  $$A_i= \{a\in A\mid a\leq x_i \trm{ but }  a\not\leq x_{i-1} \},$$
and we partially order $A$ by setting $a\lhd b$ if and only if
$a\in A_i$ and $b\in A_j$ with $i<j$.
We say $a$ is in \emph{lower level} than $b$ or $b$ is in 
\emph{higher level} than $a$ if $a\lhd b$.
Note that the level $A_i$ is an empty set if and only if $x_i$ is not
an atomic element.
A pair $(L,\De)$ is said to satisfy the \emph{level condition} if
this partial order $\unlhd$ of $A$ has the following property.
  $$\trm{If } a\lhd b_1\lhd b_2 \lhd\ldots\lhd b_k
\trm{ then } a\not\leq\Jn_{i=1}^k b_i.$$

If all elements of $\De$ are left-modular, then we say $(L,\De)$
is a \emph{left-modular} lattice.
A pair $(L,\De)$ is called an \emph{LL lattice} if it is left-modular and
satisfies the level condition.

A generalized rank function $\rho:L\ra\bbN$ is defined by
$$\rho(x)=\trm{ number of $A_i$ containing atoms less than or
equal to $x$.}$$
Note that, for any $x\in L$, we have $\rho(x)= \rho(\de(x))$
where $\de(x)$ is the maximum atomic element in $[\mh,x]$.
So $\rho(\Mh)$ is not necessary equal to $n$, the length of $\De$.

In the following we list several properties in  \cite{BS:mfl97} that we
need.
\ben
\item[(A)]
If $a$ and $b$ are distinct atoms from the same level $A_i$ in a
left-modular lattice, then $a\jn b$ is above some atom $c\in A_j$ with
$j<i$.
\item[(B)]
In an LL lattice, a set $B\sbe A$ is NBB if and only if $|B\cap A_i|\le 1$
for every $i$.
\item[(C)]
Let $B$ be an NBB set in an LL lattice.
Then every atom $a\leq\Jn B$ is in the same level as some element of $B$.
In particular, any NBB base for $x$ has exactly $\rho(x)$ atoms.
\een

Blass and Sagan generalized Stanley's Total Factorization Theorem to
LL lattices using their theory of NBB sets.
Here we present two inductive proofs for their theorem.
In the first proof we will apply Theorems~\ref{t_lmfac} as well as the
theory of NBB sets.

\bth[Blass and Sagan, \cite{BS:mfl97}] \label{t_BSlfac}
If $(L,\De)$ is an LL lattice then its characteristic polynomial factors
as
$$\chi(L,t)=\prod (t-|A_i|)$$
where the product is over all non-empty levels $A_i$.
\eth
{\bf Proof of Theorem~\ref{t_BSlfac} I.}\quad
We will induct on $n$, the length of $\De$.
The theorem is trivial when $n\le 1$.
If $A_n=\emp$, then $\rho(x_n)=\rho(x_{n-1})$ and $\mu(x)=0$ for
$x\not\leq x_{n-1}$.  Thus
$\chi(L,t)=\chi([\mh,x_{n-1}],t)$,  so we are done by induction.

If $A_n\ne\emp$, consider $b\in H(L)$.
Then, by Theorem~\ref{t_NBBmu}, $b$ must have an NBB base, say $B$.
In addition, if $b\mt x_{n-1}=\mh$ then $B\sbe A_n$ and also
$|B\cap A_n|\le 1$ by (B).
So $b=\mh$ or $b\in A_n$.
Now it suffices to check that $\tau_b$ is rank-preserving on
$H(\mh,x_{n-1})\cup\{x_{n-1}\}$ for every $b\in A_n$ since then we get
$\chi(L,t)=\chi([\mh,x_{n-1}],t)(t-|A_n|)$ by Theorem~\ref{t_lmfac}.
Because $A_n\ne\emp$ and $\rho(b)=1$, $\tau_b$ is rank-preserving on
$\{x_{n-1}\}$.
Given any $y\in H(\mh,x_{n-1})$, suppose $B$ be an NBB base for $y$.
By (B), $B'=B\cup \{b\}$ is an NBB base for $\tau_b(y)$.
Now $\rho(\tau_b(y))=|B'|=|B|+1=\rho(y)+\rho(b)$ by (C).
Hence $\rho(b,\tau_b(y))=\rho(\tau_b(y))-\rho(b)=\rho(y)=\rho(\mh,y)$.\ \Qed

In a similar way, Corollary~\ref{c_mfac} provides us with an inductive
proof for Theorem~\ref{t_mlfac}.
Note that the lattice in Theorem~\ref{t_mlfac} is graded,
so $\rho(\Mh)$ equals the length of $\De$.
Therefore the product~(\ref{e_mlfac}) is over all levels $A_i$ (including
empty ones).

We will use Theorem~\ref{t_lmgr} for the second proof.
This demonstration sidesteps the machinery of NBB sets
and reveals some properties of LL lattices in
the process.
To prepare, we need the following two lemmas.

\ble              \label{l_lmcov}
If $w$ is a left-modular element in $L$ and $v\prec w$, then
$v\jn u\preceq w\jn u$ for any $u\in L$.
\ele
\pf
Suppose not and then there exists $s\in L$ such that
$v\jn u<s<w\jn u$.
Taking the join with $w$ and using $v\jn w=w$, we get
$w\jn(v\jn u)=w\jn s=w\jn(w\jn u)$.
So we should have $w\mt(v\jn u)<w\mt s<w\mt(w\jn u)=w$ by
Theorem~\ref{t_comp}.
Combining this with $v\le w\mt(v\jn u)$, we have a contradiction to
$v\prec w$.\ \Qed

\ble       \label{l_LLforb}
If $(L,\De)$ is an LL lattice with
$\De:\mh=x_0\prec x_1\prec\ldots\prec x_n=\Mh$ and $A_n\ne\emp$,
then $([b,\Mh],\De')$ is also an LL lattice
for any $b\in A_n$ where $\De'$ consists of the distinct elements of
the multichain
$$b=x_0'\preceq x_1'\preceq x_2'\preceq\ldots\preceq x_{n-2}'\preceq
x_{n-1}'=\Mh$$
where $x_i'=x_i\jn b, 0\le i\le n-1$.
Furthermore we have $|A_i|=|A_i'|$ for such $i$, where
  $$A_i'=\{a\in A(b,\Mh)\mid a\le x_i'\trm{ but } a\not\le x_{i-1}'\}.$$
\ele
\pf
By Lemma~\ref{l_lmcov},
the chain $\De'$ is indeed saturated.
So $\De'$ is a left-modular maximal chain by Proposition~\ref{p_lmjm}.

Let $\tau(x)=\tau_b(x)=x\jn b$.
This map is surjective (see the proof of Theorem~\ref{t_lmfac}) and
order-preserving from $[\mh,x_{n-1}]$ to $[b,\Mh]$.
Also let $A=A(\mh,x_{n-1})$ and $A'=A(b,\Mh)$.
First, We prove that the map $\tau:A\ra A'$ is well-defined and bijective.
Suppose that there is an $a\in A_i$ such that $b\prec x<\tau(a)=a\jn b$ for
some $x$.
By the level condition, any atom $c\le a\jn b$ is in a level at least as
high as $a$; furthermore, if $c\in A_i$ we must have $c=a$ because of (A).
Since $x<a\jn b$ and $a\not\le x$,
any atom $d\le x$ is in a higher level than $a$.
It follows that $x_i\mt x=\mh$.
Now $b\jn(x_i\mt x)=b$ and $(b\jn x_i)\mt x\ge(b\jn a)\mt x=x$
contradicts
the left-modularity of $x_i$.
We conclude that $\tau:A\ra A'$ is well-defined.

The restriction $\tau|_{A}$ is surjective since $\tau$ is surjective and
order-preserving.
To show injectivity of $\tau|_A$, let us
suppose there are two distinct atoms $u$ and $v$ such that
$\tau(u)=\tau(v)$.
If $u$ and $v$ are from two different levels then this contradicts the
level condition.
If $u$ and $v$ are from the same level, by (A), there
exists an atom $c$ in a lower level such that
$c\le u\jn v\le\tau(u)\jn\tau(v)=\tau(u)$,
contradicting the level condition again.

Now let us prove $|A_i|=|A_i'|$.
This is trivial for $i=1$.
Let $u\in A_i$ for some nonempty $A_i$ with $2\le i\le n-1$.
It is clear that $\tau(u)\le x_i'$.
Suppose that $\tau(u)\le x_{i-1}'$, i.e., $u\jn b\le x_{i-1}\jn b$.
By the level condition, $b\jn(x_{i-1}\mt(u\jn b))=b\jn\mh=b$.  But
$(b\jn x_{i-1})\mt(u\jn b)=u\jn b>b$ contradicts the modularity of
$x_{i-1}$.
Thus $\tau(A_i)\sbe A_i'$ and then the bijectivity of $\tau|_A$
implies that $|A_i|=|A_i'|$ for all $i\le n-1$.

Since $\tau|_A$ is bijective and level-preserving,
if $\tau(a)\le\Jn_{i=1}^k \tau(b_i)$ for some
$\tau(a)\lhd \tau(b_1) \lhd \tau(b_2) \lhd\ldots\lhd \tau(b_k)$ in
$[b,\Mh]$,
then $a<a\jn b\le (\Jn_{i=1}^k b_i)\jn b$ with
$a\lhd b_1\lhd b_2\lhd\ldots\lhd b_k \lhd b$ in $L$.
Therefore $([b,\Mh],\De')$ satisfies the level condition.
\ \Qed

{\bf Proof of Theorem~\ref{t_BSlfac} II.}\quad
We will induct on $n=\ell(\De)$.
The cases $n\le 1$ and $A_n=\emp$ are handled as before.

If $A_n\neq\emp$, consider $b\in H(L)$
with $b\mt x_{n-1}=\mh$.  Then $b$ is atomic and can
only be above atoms in $A_n$.  So by (A), $b$ must be
the join of at most one atom, i.e., either $b=\mh$ or $b\in A_n$.
Thus by Lemma~\ref{l_LLforb} and induction we get, for any $b\in A_n$,
$$
\chi([b,\Mh],t)=\prod_{i\le n-1} (t-|A_i'|)
=\prod_{i\le n-1} (t-|A_i|)=\chi([\mh,x_{n-1}],t)
$$
where the product is over all non-empty $A_i$.
Applying Theorem~\ref{t_lmgr} gives
$\chi(L,t)=\chi([\mh,x_{n-1}],t)(t-|A_n|)$, so again we are done.
\ \Qed

\begin{\bib}{99}

\bibi{Bir:lt67} G. Birkhoff, ``Lattice Theory'', Third Edition,
American Mathematical Society, 1967.

\bibi{BS:mfl97} A. Blass and B. E. Sagan, M\"obius Functions of Lattices,
{\it \aim} {\bf 127} (1997), no. 1, 94--123.

\bibi{Cra:mfl66} H.~Crapo, The M\"obius function of a lattice,
{\it \jct} {\bf 1} (1966), 126--131.

\bibi{Gre:ps88} C. Greene, Posets of Shuffles,
{\it \jcta} {\bf 47} (1988), 191--206.

\bibi{Kre:pnc72}  G. Kreweras, Sur les partitions non-crois\'ees d'un cycle,
{\it \dm} {\bf 1} (1972), 333--350.

\bibi{Sta:meg71} R. P. Stanley, Modular Elements of Geometric Lattices,
{\it \alu} {\bf 1} (1971), 214--217.

\bibi{Sta:sl72} R. P. Stanley, Supersolvable Lattices,
{\it \alu} {\bf 2} (1972), 197--217.


\end{\bib}

\end{document}